\documentclass[12pt]{amsart}
\usepackage{amsmath,amssymb,amsbsy,amsfonts,latexsym,amsopn,amstext,
                                               amsxtra,euscript,amscd,bm}
\usepackage{url}
\usepackage[colorlinks,linkcolor=blue,anchorcolor=blue,citecolor=blue]{hyperref}
\usepackage{color}
\usepackage{graphics,epsfig}
\usepackage{graphicx}
\usepackage{float}

\hypersetup{breaklinks=true}

\begin{document}

\newtheorem{theorem}{Theorem}
\newtheorem{lemma}[theorem]{Lemma}
\newtheorem{claim}[theorem]{Claim}
\newtheorem{cor}[theorem]{Corollary}
\newtheorem{prop}[theorem]{Proposition}
\newtheorem{definition}{Definition}
\newtheorem{question}[theorem]{Question}
\newcommand{\hh}{{{\mathrm h}}}

\numberwithin{equation}{section}
\numberwithin{theorem}{section}
\numberwithin{table}{section}

\def\sssum{\mathop{\sum\!\sum\!\sum}}
\def\ssum{\mathop{\sum\ldots \sum}}
\def\iint{\mathop{\int\ldots \int}}

\def\squareforqed{\hbox{\rlap{$\sqcap$}$\sqcup$}}
\def\qed{\ifmmode\squareforqed\else{\unskip\nobreak\hfil
\penalty50\hskip1emrll\nobreak\hfil\squareforqed
\parfillskip=0pt\finalhyphendemerits=0\endgraf}\fi}

\newfont{\teneufm}{eufm10}
\newfont{\seveneufm}{eufm7}
\newfont{\fiveeufm}{eufm5}
%
%
\newfam\eufmfam
     \textfont\eufmfam=\teneufm
\scriptfont\eufmfam=\seveneufm
     \scriptscriptfont\eufmfam=\fiveeufm
%
%
\def\frak#1{{\fam\eufmfam\relax#1}}

\newcommand{\bflambda}{{\boldsymbol{\lambda}}}
\newcommand{\bfmu}{{\boldsymbol{\mu}}}
\newcommand{\bfxi}{{\boldsymbol{\xi}}}
\newcommand{\bfrho}{{\boldsymbol{\rho}}}

\def\fK{\mathfrak K}
\def\fT{\mathfrak{T}}

\def\fA{{\mathfrak A}}
\def\fB{{\mathfrak B}}
\def\fC{{\mathfrak C}}
\def\fM{{\mathfrak M}}

\def \balpha{\bm{\alpha}}
\def \bbeta{\bm{\beta}}
\def \bgamma{\bm{\gamma}}
\def \blambda{\bm{\lambda}}
\def \bchi{\bm{\chi}}
\def \bphi{\bm{\varphi}}
\def \bpsi{\bm{\psi}}

\def\eqref#1{(\ref{#1})}

\def\vec#1{\mathbf{#1}}


\def\cA{{\mathcal A}}
\def\cB{{\mathcal B}}
\def\cC{{\mathcal C}}
\def\cD{{\mathcal D}}
\def\cE{{\mathcal E}}
\def\cF{{\mathcal F}}
\def\cG{{\mathcal G}}
\def\cH{{\mathcal H}}
\def\cI{{\mathcal I}}
\def\cJ{{\mathcal J}}
\def\cK{{\mathcal K}}
\def\cL{{\mathcal L}}
\def\cM{{\mathcal M}}
\def\cN{{\mathcal N}}
\def\cO{{\mathcal O}}
\def\cP{{\mathcal P}}
\def\cQ{{\mathcal Q}}
\def\cR{{\mathcal R}}
\def\cS{{\mathcal S}}
\def\cT{{\mathcal T}}
\def\cU{{\mathcal U}}
\def\cV{{\mathcal V}}
\def\cW{{\mathcal W}}
\def\cX{{\mathcal X}}
\def\cY{{\mathcal Y}}
\def\cZ{{\mathcal Z}}
\newcommand{\rmod}[1]{\: \mbox{mod} \: #1}

\def\cg{{\mathcal g}}

\def\vr{\mathbf r}

\def\e{{\mathbf{\,e}}}
\def\ep{{\mathbf{\,e}}_p}
\def\eq{{\mathbf{\,e}}_q}

\def\Tr{{\mathrm{Tr}}}
\def\Nm{{\mathrm{Nm}}}

 \def\SS{{\mathbf{S}}}

\def\lcm{{\mathrm{lcm}}}

\def\({\left(}
\def\){\right)}
\def\l|{\left|}
\def\r|{\right|}
\def\fl#1{\left\lfloor#1\right\rfloor}
\def\rf#1{\left\lceil#1\right\rceil}

\def\mand{\qquad \mbox{and} \qquad}

\newcommand{\commK}[1]{\marginpar{%
\begin{color}{red}
\vskip-\baselineskip 
\raggedright\footnotesize
\itshape\hrule \smallskip B: #1\par\smallskip\hrule\end{color}}}

\newcommand{\commI}[1]{\marginpar{%
\begin{color}{magenta}
\vskip-\baselineskip 
\raggedright\footnotesize
\itshape\hrule \smallskip I: #1\par\smallskip\hrule\end{color}}}

\newcommand{\commT}[1]{\marginpar{%
\begin{color}{blue}
\vskip-\baselineskip 
\raggedright\footnotesize
\itshape\hrule \smallskip I: #1\par\smallskip\hrule\end{color}}}




\hyphenation{re-pub-lished}

\mathsurround=1pt

\def\bfdefault{b}

\def \F{{\mathbb F}}
\def \K{{\mathbb K}}
\def \Z{{\mathbb Z}}
\def \Q{{\mathbb Q}}
\def \R{{\mathbb R}}
\def \C{{\\mathbb C}}
\def\Fp{\F_p}
\def \fp{\Fp^*}

\def\Kmnp{\cK_p(m,n)}
\def\Kmnq{\cK_q(m,n)}

\def\SAMJp{\cS_p(\cA;\cM,\cJ)}
\def\SAMJq{\cS_q(\cA;\cM,\cJ)}
\def\SAJq{\cS_q(\cA;\cJ)}
\def\SAIJp{\cS_p(\cA;\cI,\cJ)}
\def\SAIJq{\cS_q(\cA;\cI,\cJ)}

\def\TWXJp{\cT_p(\cW;\varXi,\cJ)}
\def\TWXJq{\cT_q(\cW;\varXi,\cJ)}
\def\TWJq{\cT_q(\cW;\cJ)}

 \def \xbar{\overline x}

\title[Bilinear Forms with Kloosterman and Gauss Sums]{Bilinear Forms with Kloosterman and Gauss Sums}

 \author[I. E. Shparlinski] {Igor E. Shparlinski}

\address{Department of Pure Mathematics, University of New South Wales,
Sydney, NSW 2052, Australia}
\email{igor.shparlinski@unsw.edu.au}

%

\begin{abstract} We obtain several  estimates for bilinear form with Kloosterman sums. Such  results
can be interpreted as a measure of cancellations amongst with
parameters from short intervals. In particular, for certain ranges of parameters we improve
some recent results of Blomer,  Fouvry, Kowalski, Michel and Mili{\'c}evi{\'c} 
and also of Fouvry,   Kowalski and Michel. In particular, we improve the bound on the error term 
in the asymptotic formula for mixed moments of $L$-series associated 
with  Hecke eigenforms. 
 \end{abstract}


\maketitle

\section{Introduction}
\label{sec:intro}

\subsection{Background and motivation}

 Let  $q$ be  a positive 
integer. We denote the residue ring modulo $q$ by $\Z_q$ and denote the group
of units of $\Z_q$ by $\Z_q^*$. 

 For integers $m$ and $n$ we define
the Kloosterman sum
$$
\Kmnq = \sum_{x \in \Z_q^*} \eq\(mx +n\xbar \),
$$
where $\xbar$ is the multiplicative inverse of $x$ modulo $q$ and
$$
\eq(z) = \exp(2 \pi i z/q).
$$
Given a set $\cM\subseteq \Z_q^*$, an interval
$ \cJ = \{L+1, \ldots,  L+N\}\subseteq [1, q-1]$ of $N$ consecutive integers 
and a   sequence of weights $\cA = \{\alpha_m\}_{m\in \cM}$,
we define the weighted sums of Kloosterman sums
$$ 
\SAMJq = \sum_{m\in \cM} \sum_{n \in \cJ} \alpha_m \cK_q(m,n).
$$

By the Weil bound we have
$$
|\Kmnq|\le  q^{1/2+o(1)},
$$
see~\cite[Corollary~11.12]{IwKow}. Hence we immediately obtain 
 \begin{equation}
\label{eq:trivial AIJ}
\left| \SAMJq \right| \le N q^{1/2+o(1)} \sum_{m\in \cM}  |\alpha_m|.
\end{equation}
We are interested in studying cancellations amongst Kloosterman sums and
thus in improvements of the trivial bound~\eqref{eq:trivial AIJ}.

We remark that if $q=p$ is prime, then 
making the change of variable $x \mapsto nx \pmod p$, one immediately observes  that 
$\cK_p(mn,1)  = \Kmnp$, thus  we also have
$$
\SAMJp = \sum_{m\in \cM} \sum_{n \in \cJ} \alpha_m \cK_p(mn,1) .
$$
If furthermore $\cM = \cI =  \{K+1, \ldots, K+M\} \subseteq [1, p-1]$ is an interval of $M$ 
consecutive integers, we obtain the sums $\SAIJp$,  which have been studied in 
 recent works of Blomer,  Fouvry, Kowalski, Michel and Mili{\'c}e\-vi{\'c}~\cite{BFKMM1}, 
Fouvry,   Kowalski and Michel~\cite{FKM} and Shparlinski and Zhang~\cite{ShpZha}.

We have to stress that the most important part of the very deep 
work of Blomer,  Fouvry, Kowalski, Michel and 
Mili{\'c}e\-vi{\'c}~\cite{BFKMM1} is establishing the connection between bilinear forms with Kloosterman 
sums and  mixed moments of $L$-series associated 
with  Hecke eigenforms. However improving one of their ingredients allows us to improve 
one of their main results. 

The sums $\SAIJp$, including the special case of sums without weights, that is, $\cS_p\(\{1\}_{m=1}^M;\cI,\cJ\)$, as well as some other related sums,  such as 
$$
\sum_{m\in \cI} \sum_{n \in \cJ} \alpha_m \beta_n \cK_p(mn,1) \mand  
\sum_{m\in \cI}  \cK_p(m,1), 
$$
(where $\{\beta_n\}_{n\in \cJ}$ is another sequence of weights) 
appear in some applications and are also of independent interest, we refer 
to~\cite{BFKMM1,BFKMM2,FKM,FKMRRS,KMS,ShpZha,Xi-FKM} for  a wide range of 
various applications and further references.
We also recall recent results of~\cite{LSZ,Khan,WuXi} when cancellations among
Kloosterman sums are studied for moduli of special arithmetic structure.

Here we consider more general sums $\SAMJq$ extending the previously studied sums in
the following two aspects: 

\begin{itemize}
\item  the modulus  $q$ is now  an arbitrary 
positive integer;  
\item the weights are supported on   an arbitrary subset of $\cM \subseteq \Z_q^*$.
\end{itemize}
 
Using the approach of~\cite{ShpZha}, augmented with several new arguments, 
we improve and generalise previous bounds on these sums.
It is important to note, that our method does not rely on algebraic geometry 
results, as the method of~\cite{BFKMM1,FKM,KMS} and thus work modulo composite 
numbers as well as modulo primes. 

Furthermore, we use similar ideas to   the weighted sums of Gauss sums
$$
 \cG_q(\chi,n) = \sum_{x\in \Z_q^*} \chi(x) \eq(nx), 
$$ 
where $\chi$ is a  Dirichlet  character modulo $q$, see~\cite[Chapter~3]{IwKow} for a background 
on characters. More precisely, given   a   sequence of weights $\cW = \{ \omega_\chi \}_{\chi \in \varXi}$,
supported on a  subset $\varXi$ of  the set $\Omega_q^*$ of primitive Dirichlet character modulo $q$. 
we define the weighted sums of Gauss sums
$$
\TWXJq = \sum_{\chi \in \varXi} \sum_{n \in \cJ} \omega_\chi \cG_q(\chi,n).
$$
where, as before, $ \cJ = \{L+1, \ldots,  L+N\}\subseteq [1, q-1]$. 
We certainly have the full analogue of~\eqref{eq:trivial AIJ}:
$$
\left| \TWXJq \right| \le N q^{1/2+o(1)} \sum_{\chi \in \varXi}  |\omega_\chi|, 
$$
see~\cite[Equation~(3.14)]{IwKow}. The study of cancellations between Gauss sums
has been initiated by Katz and Zheng~\cite{KaZh}, in a different form as a result about the 
uniformity of distribution of their arguments. 
Furthermore, in the case of  prime $q=p$ and constant weights, one can 
obtain a nontrivial upper bound on  
$$
\cT_p\(\{1\}_{\chi \in \varXi}; \varXi, \cJ\) = \sum_{\chi \in \varXi} \sum_{n \in \cJ}  \cG_p(\chi,n)
$$
under the condition $MN \ge p^{1+\varepsilon}$ for a fixed $\varepsilon > 0$, where $M = \# \varXi$, 
 from the uniformity of distribution result of~\cite[Theorem~3]{Shp}. 
Here we obtain  stronger and more general bounds. 

As an application of our new bounds with Kloosterman sums, we  also improve the power saving in the error 
term of the asymptotic formula for mixed moments of $L$-series associated 
with  Hecke eigenforms, which refines the previous result of
Blomer,  Fouvry, Kowalski, Michel and Mili{\'c}evi{\'c}~\cite[Theorem~1.2]{BFKMM1}.

\subsection{General notation}

We  define the norms 
 $$
 \|\cA\|_\infty=\max_{m\in \cM}|\alpha_m|  \mand \|\cA\|_\sigma =\( \sum_{m\in \cM} |\alpha_m|^\sigma\)^{1/\sigma},
 $$
where  $\sigma >0$.

We always assume that  the sequence of weights $\cA = \{\alpha_m\}_{m\in \cI}$ is 
supported only on $m$ with $\gcd(m,q)=1$, that is, we have $\alpha_m = 0$  if $\gcd(m,q)>1$.

Throughout the paper,  as usual $A\ll B$  is equivalent to the inequality $|A|\le cB$
with some  constant $c>0$, which occasionally, where obvious, may depend on the real 
parameter $\varepsilon>0$ and on the integer parameter $r \ge 1$, and is absolute otherwise. 

The letter $p$ always denotes a prime number. 

\subsection{Previous results}

For 
$$
 \cI =  \{K+1, \ldots, K+M\} , \ \cJ = \{L+1, \ldots,  L+N\}\subseteq [1, q-1],
$$
the sums   $\SAIJp$
 have been estimated by Fouvry,   Kowalski and Michel~\cite[Theorem~1.17]{FKM}
 as  a part of a much more general result about sums of so-called {\it trace functions\/}.
For example, by~\cite[Theorem~1.17(2)]{FKM}, for initial intervals
 $\cI =  \{1, \ldots, M\}$ and $ \cJ = \{1, \ldots,  N\}$, we have
 \begin{equation}
\label{eqSAIJ-1}
| \SAIJp|  \le  \|\cA\|_1 p^{1+o(1)} .
 \end{equation}
 Furthermore, by a result of Blomer,  Fouvry, Kowalski, Michel 
 and Mili{\'c}evi{\'c}~\cite[Theorem~6.1]{BFKMM1},
 also for an initial interval $\cI$ and an arbitrary interval $\cJ$ with
 \begin{equation}
\label{eqSAIJ-2 Cond}
 MN \le p^{3/2} \mand M \le N^2, 
 \end{equation}
 we have
 \begin{equation}
\label{eqSAIJ-2}
\left| \SAIJp \right| \le  \( \|\cA\|_1 \|\cA\|_2\)^{1/2}  M^{1/12} N^{7/12} p^{3/4+o(1)}.
 \end{equation}

The results of~\cite{BFKMM1,FKM} are based on deep methods originating  from  
algebraic geometry, such the {\it Weil and Deligne bounds\/}, see~\cite[Chapter~11]{IwKow}. 
A much more elementary approach,  suggested in~\cite{ShpZha}, yields the estimate
\begin{equation}
\label{eqSAIJ-3}
 \SAIJp  \ll   \|\cA\|_2 N^{1/2}p.
 \end{equation}
In particular, we see that the approach of~\cite{ShpZha} improves the 
bounds from~\cite{BFKMM1,FKM}
for 
$$
N < M p^{-\varepsilon} \mand M^4 N \ge p^{3 + \varepsilon}
$$
with any fixed $\varepsilon>0$. 
 
For the purpose of comparison between previous results and our new bounds, we
rewriting the bounds~\eqref{eq:trivial AIJ}, \eqref{eqSAIJ-1}, \eqref{eqSAIJ-2} and~\eqref{eqSAIJ-3}
in terms of $\|\cA\|_\infty$ and combine them in one bound
 \begin{equation}
 \label{eq:BoundComb}
\begin{split}
\SAIJp \ll   \|\cA\|_\infty \min\{MN p^{1/2}, M p &,M^{5/6} N^{7/12} p^{3/4},\\
&  \qquad M^{1/2} N^{1/2} p\}p^{o(1)}, 
\end{split}
\end{equation}
where we also ignore the necessary condition~\eqref{eqSAIJ-2 Cond} for the bound~ \eqref{eqSAIJ-2}
to apply. 

For sums of Gauss sums, no general results have been known. However, for a prime $q=p$ as we have mentioned, 
one can derive a nontrivial bound in the case of constant weights from~\cite[Theorem~3]{Shp}
and in fact for more general sums  with the summation over $n$ over an arbitrary set $\cN \subseteq \Z_p$.

\section{New results}

\subsection{Bounds of sums of Kloosterman sum}

We remark that our bounds only involve the norms of the weights $\cA$ 
but do not explicitly depend on the size of the set $\cM$ on which they are supported. 
Hence, without loss of generality, we can assume that $\cM = \Z_q^*$ 
and thus we simplify the notation as 
$$
\SAJq = \sum_{m\in  \Z_q^*} \sum_{n \in \cJ} \alpha_m \cK_q(m,n).
$$

\begin{theorem}
\label{thm:SAJq} For any integer $q \ge 1$,  we have,
$$
 \SAJq  \ll \( \|\cA\|_1  \|\cA\|_2\)^{1/2} \(N^{1/8}q   +N^{1/2}q^{3/4}\)q^{o(1)}. 
$$
\end{theorem}

\begin{figure}[H]
  \centering
  \includegraphics[scale=0.6]{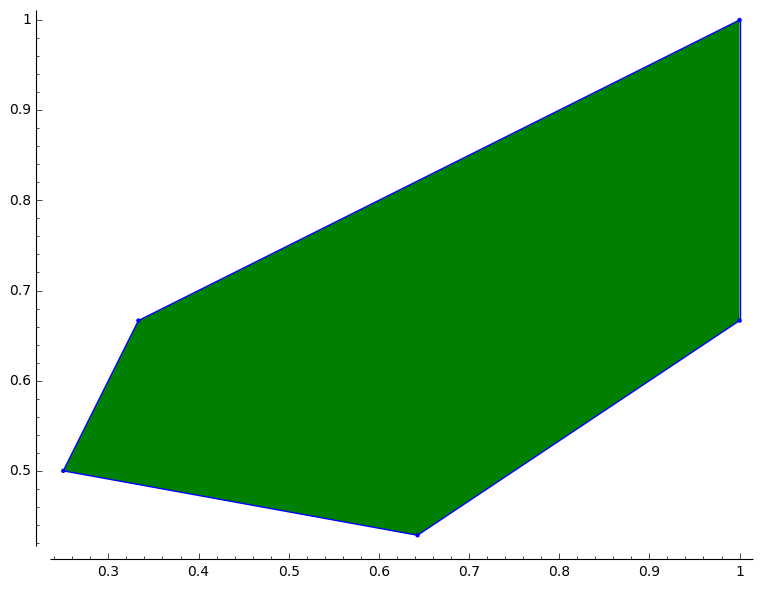}
  \caption{Polygon of relative sizes of $M$ and $N$ defined by~\eqref{eq:Improve}.}
  \label{fig:Improve}
\end{figure}

Returning to our original settings and assuming that the sequence of weights $\cA$ is 
supported on a set $\cM \subseteq \Z_q^*$ of size $M$, 
we can rewrite that bound of Theorem~\ref{thm:SAJq} 
in terms of $\|\cA\|_\infty$  as
 \begin{equation}
 \label{eq:SimpleBound}
 \SAMJq  \ll \|\cA\|_\infty M^{3/4}  \(N^{1/8}  q  +N^{1/2}q^{3/4}\) q^{o(1)}.
\end{equation}
We see that the bound~\eqref{eq:SimpleBound}, besides being more general than~\eqref{eq:BoundComb}, 
also gives a better result, provided that
 \begin{equation}
 \label{eq:Improve}
M^2 N^7 \ge q^4, \quad M^2 N^{11} \ge q^6, 
 \quad qM \ge N^2, \quad N^3 \ge M^2 \ge N.
\end{equation}
Writing $M = q^\mu$ and $N=q^\nu$ we see that the conditions~\eqref{eq:Improve}
define a polygon with vertices 
$$
(1/4, 1/2),  \quad(1/3, 2/3), \quad (1, 1), \quad (1, 2/3), \quad (9/14, 3/7), 
$$
in the $(\mu,\nu)$-plane, 
see also Figure~\ref{fig:Improve}. The most important for applications fact is that the point $(1/2, 1/2)$ is
an interior point of this polygon. On the other hand, we  also remark that although we have ignored the 
restriction~\eqref{eqSAIJ-2 Cond} for the bound~\eqref{eqSAIJ-2} to hold, taking it into account does not 
increase the region where~\eqref{eq:SimpleBound} improves the previously known bounds.

In particular, even ignoring the differences in the generality,  
in the situation  where all previous bounds apply, 
when $M \sim N \sim p^{1/2}$ (which is a critical range for several applications), 
 Theorem~\ref{thm:SAJq}
implies the bound $p^{3/2 -1/16 +o(1)}$ instead of $p^{3/2 -1/24 +o(1)}$ given by~\eqref{eq:BoundComb}
(which actually comes  from~\eqref{eqSAIJ-2}). 

We also obtain a stronger bound  on average over $q$ in a dyadic interval $[Q,2Q]$:

\begin{theorem}
\label{thm:SAJq Aver}   For any fixed real $\varepsilon > 0$ and integer $r\ge 2$, for any sufficiently large  $Q \ge 1$, for all but  at most $Q^{1-2r \varepsilon + o(1)}$ integers 
$q\in [Q,2Q]$ we have 
$$
 \SAJq  \ll   \|\cA\|_1^{1-1/r}  \|\cA\|_2^{1/r}  
 \(q   +N^{1/2}q^{1/2+1/2r}\)q^{\varepsilon+o(1)}.
$$
\end{theorem}

Again, when the sequence of weights $\cA$ is 
supported on a set $\cM \subseteq \Z_q^*$ of size $M$, 
we  can rewrite that bound of Theorem~\ref{thm:SAJq Aver} 
in terms of $\|\cA\|_\infty$  as
$$
 \SAMJq  \ll \|\cA\|_\infty M^{1-1/2r } \(q   +N^{1/2}q^{1/2+1/2r}\)q^{\varepsilon+o(1)}.
$$

\subsection{Bounds of sums of Gauss sum}

For the sums of Gauss sums we also have bounds  
that only involve the norms of the weights $\cW$ 
but do not explicitly depend on the size of the set $\varXi$ on which they are supported. 
Hence, without loss of generality, we can assume that $\varXi = \Omega_q^*$ 
and thus we simplify the notation as 
$$
\TWJq = \sum_{\chi \in \Omega_q^*} \sum_{n \in \cJ} \omega_\chi \cG_q(\chi,n).
$$

\begin{theorem}
\label{thm:TWJq} For any integer $q \ge 1$,  we have,
$$
\TWJq \ll \( \|\cW\|_1  \|\cW\|_2\)^{1/2} \(q   +N^{1/2}q^{3/4}\)q^{o(1)}. 
$$
\end{theorem}

Finally, we also have
\begin{theorem}
\label{thm:TWJq Aver}   For any fixed real $\varepsilon > 0$ and integer $r\ge 2$, for any sufficiently large  $Q \ge 1$, for all but at most  $Q^{1-2r \varepsilon + o(1)}$ integers 
$q\in [Q,2Q]$ we have 
$$
\TWJq   \ll   \|\cW\|_1^{1-1/r}  \|\cW\|_2^{1/r}  
 \(q   +N^{1/2}q^{1/2+1/2r}\)q^{\varepsilon+o(1)}.
$$
\end{theorem}

If  the sequence of weights $\cW$ is 
supported on a set $\varXi \subseteq \Omega_q^*$ of size $M$, 
then we have 
$$
\TWXJq  \ll \|\cW\|_\infty M^{1-1/2r } \(q   +N^{1/2}q^{1/2+1/2r}\)q^{\varepsilon+o(1)}.
$$
where we can take $r = 2$ under the conditions of Theorem~\ref{thm:TWJq}
and any  integer $r\ge 2$ under the conditions of Theorem~\ref{thm:TWJq Aver}.

%
%

\subsection{Applications}

Inserting the bound of Theorem~\ref{thm:SAJq}  (in the special case of a  prime $q =p$ and 
the weights supported on the initial interval 
$\cI_0 =  \{1, \ldots, M\}$) in the argument 
of the proof of~\cite[Theorem~1.2]{BFKMM1},  we  improve the error 
term in the asymptotic formula for mixed moments of $L$-series associated 
with  Hecke eigenforms. 

To be more precise we need to recall some definitions from~\cite{BFKMM1}. 

Let $f$ be a Hecke eigenform and let $\chi$ be a primitive Dirichlet character modulo $q$. 
We define the following $L$-series by the absolutely converging for $\Re s > 1$ series:
$$
L(f\otimes \chi, s) =\sum_{n=1}^\infty \frac{\lambda_f(n)\chi(n)}{n^s} \quad \text{and} \quad 
L(f, s) =\sum_{n=1}^\infty \frac{\lambda_f(n)}{n^s}, \qquad \Re s > 1,
$$
where $\{\lambda_f(n)\}_{n=1}^\infty$ are the  Hecke eigenvalues of $f$. We also define
the partial derivative at $s=1/2$
$$
F(z) = \frac{\partial}{\partial s}\Big\vert_{s=1/2} E(z;s)
$$
of the   {\it Eisenstein series\/} defined by
$$
E(z;s) = \frac{1}{2}\sum_{\substack{c, d \in \Z\\\gcd(c,d)=1}} \frac{(\Im z)^s}{|cz+d|^{2s}}, 
\qquad \Re s > 1, \ \Im z > 0, 
$$
Finally, we define the average value 
$$
\fM_{f,F}(q) = \frac{1}{\#  \Omega_q^*} \sum_{\chi \in \Omega_q^*} L(f\otimes \chi, 1/2) \overline {L(F\otimes \chi, 1/2)}, 
$$
where, as before,  $\Omega_q^*$ denotes the set of all primitive Dirichlet characters modulo $q$. 
Clearly, $\# \Omega_p^* = p-2$ for a prime $p$. As usual, $\zeta(s)$ denotes the Rieman 
zeta-function. 

\begin{theorem}
\label{thm:MixMom}   For any prime $p$,  we have,
$$
\fM_{f,F}(p) = \frac{L(f,1)^2}{\zeta(2)} + O\(p^{-1/64 + o(1)}\),
$$
where the implied constant depends only on $f$. 
\end{theorem}

Theorem~\ref{thm:MixMom} improves the error term
$p^{-1/68 + o(1)}$ of~\cite[Theorem~1.2]{BFKMM1}.

 \section{Preliminaries} 
 
  \subsection{Equations and congruences with reciprocals} 
 \label{sec:CongEqReipr}  
 An important tool in our argument is an upper bound on the number 
of solutions   $J_r(q;K)$ to the  congruence   
$$\frac{1}{x_1}+ \ldots+ \frac{1}{x_r}\equiv \frac{1}{x_{r+1}}+ \ldots+\frac{1}{x_{2r}}\pmod q, \quad 
1 \le x_1,  \ldots, x_{2r} \le K. 
$$
where $r =1, 2, \ldots$. 

For arbitrary $q$ and $K$, good upper bounds on  $J_r(q;K)$  
are known only for $r =2$ (and of course in the trivial case $r =1$) and are due to 
Heath-Brown~\cite[Page~368]{H-B1}
(see the bound on the sums of quantities $m(s)^2$ in the 
notation of~\cite{H-B1}).  More precisely, we have:

\begin{lemma}
\label{lem:Cayley} For $1 \le K \le q$ we have 
$$
J_2(q;K) \le \(K^{7/2} q^{-1/2} + K^2\)q^{o(1)}.
$$
\end{lemma}

It is also shown by Fouvry and  Shparlinski~\cite[Lemma~2.3]{FouShp} that the bound of 
Lemma~\ref{lem:Cayley} 
can be improved on average over $q$ in a dyadic interval $[Q,2Q]$. The same argument 
also works for $J_r(q;K)$ without any changes. 

Indeed,  let  $J_r(K)$ be the number 
of solutions   to  the  equation 
$$
\frac{1}{x_1}+ \ldots+ \frac{1}{x_r} = \frac{1}{x_{r+1}}+ \ldots+\frac{1}{x_{2r}}, \qquad 
1 \le x_1,  \ldots, x_{2r} \le K, 
$$
where $r =1, 2, \ldots$. We recall  that by the result of Karatsuba~\cite{Kar} (presented in 
the proof of~\cite[Theorem~1]{Kar}), see also~\cite[Lemma~4]{BouGar1}, we have:

\begin{lemma}
\label{lem:Recipr Eq} For any fixed positive integer $r$,  we have 
$$
J_r(K) \le K^{r + o(1)}.
$$
\end{lemma}

Now repeating the argument of the proof of~\cite[Lemma~2.3]{FouShp} and using 
Lemma~\ref{lem:Recipr Eq}  in the appropriate place,  we obtain:

\begin{lemma}
\label{lem:Cayley Aver} For any   fixed positive integer $r$ and sufficiently large integers $1 \le K \le Q$, 
we have 
$$
\frac{1}{Q}\sum_{Q \le q \le 2Q} J_{r,q}(K) \le \(K^{2r} Q^{-1} + K^r\)Q^{o(1)}.
$$
\end{lemma}

%
%

  \subsection{Equations and congruences with products} 
  
Let  $R_r(q;K)$ be the number of  solutions to  the  congruence   
$$
x_1  \ldots x_r \equiv x_{r+1}\ldots   x_{2r}\pmod q, \quad 
1 \le x_1,  \ldots, x_{2r} \le K. 
$$
where $r =1, 2, \ldots$. 
The proof of Theorem~\ref{thm:TWJq}  uses 
the bound of Friedlander and Iwaniec~\cite[Lemma~3]{FrIw} on $R_2(q;K)$ 
(which is formulated as a bound on  the 4th moment of character sums), 
see also~\cite[Theorem~2]{CochShi}. We present it in a simplified form.

\begin{lemma}
\label{lem:ProdCong} For $1 \le K \le q$ we have 
$$
R_2(q;K) \le \(K^4q^{-1} + K^2\)q^{o(1)}.
$$
\end{lemma}

Now  let  $R_r(K)$ be the number 
of solutions  to  the  equation 
$$
x_1  \ldots x_r = x_{r+1}\ldots  x_{2r}, \qquad 
1 \le x_1,  \ldots, x_{2r} \le K, 
$$
where $r =1, 2, \ldots$. We classical bound on the divisor function immediately 
implies:

\begin{lemma}
\label{lem:Prod Eq} For any fixed positive integer $r$,  we have 
$$
R_r(K) \le K^{r + o(1)}.
$$
\end{lemma}

Thus using Lemma~\ref{lem:Prod Eq} instead of Lemma~\ref{lem:Recipr Eq}  we obtain an 
analogue of Lemma~\ref{lem:Cayley Aver}:

\begin{lemma}
\label{lem:ProdCong Aver} For any   fixed positive integer $r$ and sufficiently large integers $1 \le K \le Q$, 
we have 
$$
\frac{1}{Q}\sum_{Q \le q \le 2Q} R_{r,q}(K) \le \(K^{2r} Q^{-1} + K^r\)Q^{o(1)}.
$$
\end{lemma}

\section{Proofs}

\subsection{Proof of Theorem~\ref{thm:SAJq}}


For an integer $u$ we define
$$
\langle u\rangle_q = \min_{k \in \Z} |u - kq|
$$
as the distance to the closest integer, 
 which is a multiple of $q$.

Changing the order of summation and then changing the variable
$x \mapsto \xbar$, we obtain
 \begin{equation*}
\begin{split}
 \SAJq &  = \sum_{x=1}^{p-1} \sum_{m\in \Z_q^*} 
 \alpha_m \eq(mx) \sum_{n \in \cJ} \eq(n\xbar)\\
& =\sum_{x\in \Z_q^*} \sum_{m\in  \Z_q^*} \alpha_m \eq(m\xbar) \sum_{n \in \cJ} \eq(n x).
\end{split}
\end{equation*}
Hence
$$
 \SAJq   = \sum_{m\in  \Z_q^*} \sum_{x\in \Z_q^*} \alpha_m
 \gamma_x  \eq(m \xbar),
$$
where
$$
|\gamma_x| \le  \min\left\{N,  \frac{q}{\langle u\rangle_q}\right \}.
$$

We now set  $I = \rf{\log (N/2)}$ and define $2(I+1)$ the sets
 \begin{equation*}
\begin{split}
\cX_0^{\pm} &   =\{x \in \Z~: ~ 0 < \pm x \le q/N\}, \\
\cX_i^{\pm}  & = \{x \in \Z~: ~\min\{q/2, e^{i} q/N\}\ge \pm x > e^{i-1}q/N \}, \quad i = 1, \ldots, I.
\end{split}
\end{equation*}
Therefore, 
 \begin{equation}
\label{eq:BAIJ S0i}
 \SAJq \ll   \sum_{i=0}^I \(|S_i^+| + |S_i^-|\),
\end{equation}
where
$$
S_i ^{\pm} =  \sum_{m\in  \Z_q^*} \sum_{x \in \cX_i^\pm} \alpha_m \gamma_x  \eq(m \xbar), \qquad i = 0, \ldots, I.
$$

Blow we present the argument in a general form with an arbitrary $r \ge 2$. 
We then apply it with $r =2$ since we use Lemma~\ref{lem:Cayley}. However 
in the proof of Theorem~\ref{thm:SAJq Aver} we use it in full generality.

Le us fix some integer $r \ge 2$. Writing 
$$
|S_i^\pm|  \le   \sum_{m\in  \Z_q^*} \left|  \alpha_m\right|^{(r-1)/r}
 \left|  \alpha_m^2\right|^{1/2 r }
\left |\sum_{x \in \cX_i^\pm} \alpha_m \gamma_x  \eq(m \xbar)\right|, 
$$
by the H{\"o}lder inequality, for every $i =0, \ldots, I$ and every choice of the 
sign `$+$' or `$-$',  we obtain
 \begin{equation}
 \label{eq:Holder S}
\begin{split}
|S_i^\pm|   & \le  
\( \sum_{m\in  \Z_q^*} |\alpha_m|\)^{1 -1/r}
\( \sum_{m\in  \Z_q^*} |\alpha_m|^2\)^{1/2r} \\
& \qquad \qquad \qquad \qquad \qquad \(  \sum_{m\in  \Z_q^*}\left| \sum_{x \in \cX_i^\pm} 
 \gamma_x  \eq(m \xbar)\right|^{2r} \)^{1/2r}\\
&  = \|\cA\|_1^{1-1/r}  \|\cA\|_2^{1/r} \( \sum_{m\in  \Z_q^*}\left| \sum_{x \in \cX_i^\pm} 
 \gamma_x  \eq(m \xbar)\right|^{2r} \)^{1/2r}. 
\end{split}
\end{equation}
Extending the summation over $m$ to the whole ring $\Z_q$, 
opening up the inner sum, changing the order of summation and 
using the orthogonality of exponential functions, we obtain 
 \begin{equation*}
\begin{split}
\sum_{m\in  \Z_q^*}&\left| \sum_{x \in \cX_i^\pm} 
 \gamma_x  \eq(m \xbar)\right|^{2r} \\
 &   \le \sum_{m\in \Z_q}  \ssum_{x_1,  \ldots,  x_{2r} \in \cX_i^\pm} 
\prod_{j=1}^r  \gamma_{x_j}  \overline{\gamma_{x_{r + j}} }
 \eq\(m \sum_{j=1}^r \(\xbar_j - \xbar_{r + j}\)\)\\
 &   \le  \ssum_{x_1,  \ldots,  x_{2r} \in \cX_i^\pm} 
\prod_{j=1}^r  \gamma_{x_j}  \overline{\gamma_{x_{r + j}} }
\sum_{m\in \Z_q}  \eq\(m \sum_{j=1}^r \(\xbar_j - \xbar_{r + j}\)\)\\
 & =  q  \ssum_{\substack{x_1,  \ldots,  x_{2r} \in \cX_i^\pm\\ 
 \xbar_1 + \ldots + \xbar_r \equiv \xbar_{r+1}+ \ldots + \xbar_{2r} \pmod q}}
\prod_{j=1}^r  \gamma_{x_j}  \overline{\gamma_{x_{r + j}} }.
\end{split}
\end{equation*}

We also observe that for $x \in \cX_i^\pm$ we have 
$$
|\gamma_x| \ll e^{-i} N. 
$$
Hence,  for every $i=0, \ldots, I$ we have
\begin{equation}
\label{eq:Si J}
\begin{split}
\sum_{m\in  \Z_q^*}&\left| \sum_{x \in \cX_i^\pm} 
 \gamma_x  \eq(m \xbar)\right|^{2r}\\
  &\qquad  \ll   e^{-2r i}  N^{2r} q 
 \ssum_{\substack{x_1,  \ldots,  x_{2r} \in \cX_i^\pm\\ 
 \xbar_1 + \ldots + \xbar_r \equiv \xbar_{r+1}+ \ldots + \xbar_{2r} \pmod q}}1\\
& \qquad \le    e^{-2r i}  N^{2r} q  J_r(q;\fl{e^{i}q/N}).
 \end{split}
\end{equation}

 Now using~\eqref{eq:Si J} with $r = 2$ and revoking Lemma~\ref{lem:Cayley} we obtain 
  $$
 \sum_{m\in  \Z_q^*}\left| \sum_{x \in \cX_0^\pm} 
 \gamma_x  \eq(m \xbar)\right|^{4}   \le    e^{-4i} N^{4}   \(N^{-7/2}q^3   + N^{-2}q^2\) q, 
 $$
 Nw we see from~\eqref{eq:Holder S} that 
 \begin{equation}
\label{eq:Si}
\begin{split}
|S_i^{\pm}| & \le  \( \|\cA\|_1  \|\cA\|_2\)^{1/2}\\
& \qquad \qquad   e^{-i}  N   \(e^{7i/2}N^{-7/2}q^3 + e^{2i}  N^{-2}q^2  \)^{1/4}q^{1/4+o(1)} \\
& \le e^{-i/8}  \( \|\cA\|_1  \|\cA\|_2\)^{1/2} \(N^{1/8}q   +N^{1/2}q^{3/4}\)q^{o(1)}. 
\end{split}
\end{equation}
Therefore,
\begin{equation}
\label{eq: sum Si}
\sum_{i=0}^I |S_i^\pm | \le  \( \|\cA\|_1  \|\cA\|_2\)^{1/2} \(N^{1/8}q   +N^{1/2}q^{3/4}\)q^{o(1)}. 
\end{equation}
Substituting~\eqref{eq: sum Si} in~\eqref{eq:BAIJ S0i}, 
we obtain the result.

\subsection{Proof of Theorem~\ref{thm:SAJq Aver}}

We proceed as in the proof of Theorem~\ref{thm:SAJq}, in particular, 
we set $I = \rf{\log (N/2)}$. We also define  $K_i = \fl{2 e^{i}Q/N}$ 
and replace  $J_r(q;\fl{e^{i}q/N})$ 
with $J_r(q;K_i)$ in~\eqref{eq:Si J},  $i=0, \ldots, I$.

We know see that by Lemma~\ref{lem:Cayley Aver} for every $i=0, \ldots, I$ for all
but at most $Q^{1-2r \varepsilon + o(1)}$ integers  $q \in [Q,2Q]$ we have 
\begin{equation}
\label{eq: Bound JqKi}
J_{r, q}(K_i) \le \(K_i^{2r} q^{-1} + K_i^r\)Q^{2 r \varepsilon}.
\end{equation}
Since $I = Q^{o(1)}$, we see that for all
but at most $Q^{1-2r \varepsilon + o(1)}$ integers  $q \in [Q,2Q]$, 
the bound~\eqref{eq: Bound JqKi} holds for all $i=0, \ldots, I$ 
simultaneously.  For every such $q$,  using~\eqref{eq: Bound JqKi}, instead of the 
bound of  Lemma~\ref{lem:Cayley},  we obtain 
\begin{equation*}
\begin{split}
|S_i^{\pm}| & \ll  \|\cA\|_1^{1-1/r}  \|\cA\|_2^{1/r} e^{-i}  N   \(e^{2r i} N^{-2r} q^{2r-1}   + e^{r i} N^{-r}q^r\)^{1/2r} q^{1/2 r}Q^{\varepsilon} \\
& \ll  \|\cA\|_1^{1-1/r}  \|\cA\|_2^{1/r}  \(q   +N^{1/2}q^{1/2+1/2r}\)Q^{\varepsilon}
\end{split}
\end{equation*}
instead of~\eqref{eq:Si} for every $i=0, \ldots, I$.  Since $I = Q^{o(1)}$, the result now follows.

\subsection{Proof of Theorem~\ref{thm:TWJq}}

We define  $I$, the sets $\cX_i^{\pm}$, $i = 0, \ldots, I$,  
and the quantities $\gamma_x$, $x \in \Z_q^*$,  as in the proof of Theorem~\ref{thm:SAJq}. 
Therefore, instead of~\eqref{eq:BAIJ S0i} we have
 \begin{equation}
\label{eq:BWIJ T0i}
 \TWJq \ll   \sum_{i=0}^I \(|T_i^+| + |T_i^-|\),
\end{equation}
where
$$
T_i ^{\pm} =  \sum_{\chi \in \Omega_q^*} \sum_{x \in \cX_i^\pm} \omega_\chi  \gamma_x  \chi(x), \qquad i = 0, \ldots, I.
$$
Furthermore, we have the following analogue of~\eqref{eq:Holder S}
 \begin{equation}
\label{eq:Holder T}
|T_i^\pm|    \le    \|\cW\|_1^{1-1/r}  \|\cW\|_2^{1/r} \( \sum_{\chi \in \Omega_q^*} \left| \sum_{x \in \cX_i^\pm} 
 \gamma_x  \chi(x)\right|^{2r} \)^{1/2r}, 
\end{equation}
from which, by the  orthogonality of characters,  we derive an analogue of~\eqref{eq:Si J}. 
More precisely,   for every $i=0, \ldots, I$ we obtain
\begin{equation}
\label{eq:Ti R}
\sum_{\chi \in \Omega_q^*} \left| \sum_{x \in \cX_i^\pm} 
 \gamma_x  \chi(x)\right|^{2r}  \le    e^{-2r i}  N^{2r} q  R_r(q;\fl{e^{i}q/N}).
\end{equation}

 Now using~\eqref{eq:Ti R} with $r = 2$ and revoking Lemma~\ref{lem:ProdCong} we obtain 
  $$
 \sum_{m\in  \Z_q^*}\left| \sum_{x \in \cX_0^\pm} 
 \gamma_x  \eq(m \xbar)\right|^{4}   \le    e^{-4i} N^{4}   \(N^{-4}q^3   + N^{-2}q^2\) q, 
 $$
 and thus by~\eqref{eq:Holder T}, we have 
 \begin{equation*}
\begin{split}
|T_i^{\pm}| & \le  \( \|\cW\|_1  \|\cW\|_2\)^{1/2}\\
& \qquad \qquad   e^{-i}  N   \(e^{4i}N^{-4}q^3 + e^{2i}  N^{-2}q^2  \)^{1/4}q^{1/4+o(1)} \\
& \le   \( \|\cW\|_1  \|\cW\|_2\)^{1/2} \(q   +N^{1/2}q^{3/4}\)q^{o(1)}. 
\end{split}
\end{equation*}
Therefore,
\begin{equation}
\label{eq: sum Ti}
\sum_{i=0}^I |T_i^\pm | \le  \( \|\cW\|_1  \|\cW\|_2\)^{1/2} \(q   +N^{1/2}q^{3/4}\)q^{o(1)}. 
\end{equation}
Substituting~\eqref{eq: sum Ti} in~\eqref{eq:BWIJ T0i}, 
we obtain the result.

\subsection{Proof of Theorem~\ref{thm:TWJq Aver}}

We proceed as in the proof of Theorem~\ref{thm:TWJq}, 
using Lemma~\ref{lem:ProdCong Aver} instead of Lemma~\ref{lem:ProdCong} 
in the appropriate place (see also the proof of Theorem~\ref{thm:SAJq Aver}).

\subsection{Proof of Theorem~\ref{thm:MixMom}}

We simply incorporate the bound of  Theorem~\ref{thm:SAJq} in the arsenal of 
bounds used in~\cite[Section~6.4.1]{BFKMM1}. Augmenting the {\tt Mathematica}
code, provided in~\cite[Section~7.4.1]{BFKMM1},  with this new bound, we see that 
the contribution from the  terms considered 
 in~\cite[Section~6.4.1]{BFKMM1} can be estimated as  $p^{-1/52 + o(1)}$.
 Hence, now the error term is dominated by the terms treated 
 in~\cite[Section~6.4.2]{BFKMM1},  which contribute  at most $p^{-1/64 + o(1)}$. 
 The result now follows.

\section{Comments and further applications}

We note that the error term of Theorem~\ref{thm:MixMom} is now dominated by the terms
whose treatment is free of any use of sums of Kloosterman sums, 
see~\cite[Section~6.4.2]{BFKMM1}. Hence no further improvement is possible until 
this part is refined. For example,  at the moment we cannot  take any advantage 
of Theorem~\ref{thm:SAJq Aver} in this context. 

We also recall that for a prime $q=p$ and small $K$, a series of bounds on $J_r(p,K)$ 
have been given by Bourgain and Garaev~\cite{BouGar1,BouGar2}. These bounds can also be 
used in the argument of the proof of Theorem~\ref{thm:SAJq}, leading to a series 
of estimates when $N$ is close to $q$. Similarly, the bounds of~\cite{BGKS1,BGKS2}
on  $R_r(p,K)$
can  be used in the argument of the proof  of Theorem~\ref{thm:TWJq}.

We note the suggested here approach can be applied to  many other families of bilinear 
sums of the form
$$
\cS_{k,q}(\cA;\cJ)=  \sum_{m\in  \Z_q^*} \sum_{n \in \cJ} \alpha_m  \sum_{x \in \Z_q^*}  \eq(mx^{-k}+ nx)
$$
for an integer $k\ge 1$ (generalising the sums $\cS_{1,q}(\cA;\cJ) = \SAJq$). 
 Indeed,  instead of Lemma~\ref{lem:Cayley},  in the appropriate place of the argument, 
one simply uses the bound of Heath-Brown~\cite[Lemma~1]{H-B2}
for $k=2$ and a more general bound of Pierce~\cite[Theorem~4]{Pierce} 
for arbitrary integer $k\ge 1$, see also~\cite[Proposition~1]{BouGar1}. 
We note  that the sums $\cS_{2,q}(\cA;\cJ)$ and related sums
have been estimated by Nunes~\cite[Theorems~1.2 and~1.3]{Nun} as a tool
in investigating the distribution of squarefree numbers in arithmetic progressions. 
Nunes~\cite{Nun} uses the method of~\cite{BFKMM1,KMS} and it is very plausible that 
the method of this work, and also the method of~\cite{ShpZha} in the case of constant weights, may lead
to stronger results. 

Finally, it is easy to see that one can also use the same approach to estimate the sums 
$$
\sum_{m\in  \Z_p^*} \sum_{n \in \cJ} \alpha_m  \sum_{x \in \Z_t}  \ep(mg^x) \e_t(nx), 
$$
where $g$ is an integer of multiplicative order $t$ modulo a prime $p$ and 
$ \cJ = \{L+1, \ldots,  L+N\}\subseteq [1, t]$. In this case, instead of Lemma~\ref{lem:Cayley},
one uses the bound 
$$
\#\left\{1 \le u,v,x,y \le K~:~ g^u+g^v \equiv g^x + g^y\pmod p\right\} \ll K^{5/2}, 
$$
for any $K \le t$, which follows immediately from a result 
of Roche-Newton,  Rudnev and  Shkredov~\cite[Theorem~6]{RNRS}
(see also the proof of~\cite[Theorem~18]{RNRS}). 

\section*{Acknowledgement}

The authors would like to thank Ramon Nunes  for  useful discussions and 
in particular for the information about the results of Heath-Brown~\cite[Lemma~1]{H-B2}
and Pierce~\cite[Theorem~4]{Pierce}. 

%
%
This work was  supported   by ARC Grant~DP140100118. 

\end{document}